\documentclass[twocolumn,amsmath,amssymb,aps,secnumarabic,%
    nofootinbib,superscriptaddress,floatfix]{revtex4}
\usepackage{times,mathptmx}

\renewcommand{\mspace}[1]{\mskip #1mu plus 1mu minus 1mu}
\newcommand{\deq}{\mspace{5} = \mspace{5}}
\newcommand{\drel}[1]{\mspace{5} #1 \mspace{5}}

\newcommand{\opshrink}[2]{\mskip -#1mu minus -#1mu #2 \mskip -#1mu minus -#1mu}
\newcommand{\squish}[1]{\opshrink{3}{#1\nobreak}}
\newcommand{\relshrink}[2]{\mskip -#1mu plus -2mu #2 \mskip -#1mu plus -2mu}
\newcommand{\relshrnk}[1]{\relshrink{3}{#1\nobreak}}
\newcommand{\bcdisplaydiv}{\mspace{-5} \Bigg/ \mspace{-5}}
\newcommand{\PG}{\mathop{\rm PG}}
\newcommand{\AG}{\mathop{\rm AG}}
\newcommand{\GF}{\mathop{\rm GF}}
\newcommand{\qbinomial}[2]{\mathchoice{{{#1}\brack {#2}}_q}%
    {\smash{{{#1}\brack {#2}}_q}}{{{#1}\brack {#2}}_q}{{{#1}\brack {#2}}_q}}

\newtheorem{theorem}{Theorem}

\begin{document}
\title{New constructions for covering designs}
\newcommand{\CCR}
    {Center for Communications Research, 4320 Westerra Ct., San Diego, CA 92121}
\author{Daniel M. Gordon}
\affiliation{\CCR}
\author{Greg Kuperberg}
\affiliation{Department of Mathematics, University of California, One Shields
    Avenue, Davis, CA 95616}
\author{Oren Patashnik}
\affiliation{\CCR}

\date{16 February 1995}

\begin{abstract}
A $(v,k,t)$ {\em covering design}, or {\em covering}, is a family of
$k$-subsets, called blocks, chosen from a $v$-set, such that each $t$-subset is
contained in at least one of the blocks. The number of blocks is the covering's
{\em size}, and the minimum size of such a covering is denoted by $C(v,k,t)$.
This paper gives three new methods for constructing good coverings: a greedy
algorithm similar to Conway and Sloane's algorithm for lexicographic
codes~\cite{lex}, and two methods that synthesize new coverings from
preexisting ones. Using these new methods, together with results in the
literature, we build tables of upper bounds on $C(v,k,t)$ for $v \leq 32$,\, $k
\leq 16$, and $t \leq 8$.%
\end{abstract}
\maketitle

\section{Introduction}

Let the covering number $C(v,k,t)$ denote the smallest number of $k$-subsets
of a $v$-set that cover all $t$-subsets. These numbers have been studied
extensively. Mills and Mullin~\cite{mills-mullin} give known results and many
references. Hundreds of papers have been written for particular values of
$v$,~$k$, and~$t$. The best general lower bound on $C(v,k,t)$, due to
Sch\"onheim~\cite{schonheim}, comes from the following inequality:
\begin{theorem}
\label{schonheim-one-level}
$$C(v,k,t) \drel{\geq} \Bigl\lceil \frac{v}{k} \,
C(v\squish{-}1, k\squish{-}1, t\squish{-}1) \Bigr\rceil \,.$$
\end{theorem}
Iterating this gives the Sch{\"o}nheim bound $C(v,k,t) \geq L(v,k,t)$, where
$$ L(v,k,t) \deq \Bigl\lceil
\frac{v}{k} \Bigl\lceil \frac{v-1}{k-1} \ldots
\Bigl\lceil \frac{v-t+1}{k-t+1} \Bigl\rceil \ldots \Bigl\rceil \Bigl\rceil \,.
$$
Sometimes a lower bound of de~Caen~\cite{decaen-lower-bound} is slightly
better than the Sch{\"o}nheim bound when $k$~and~$t$ are not too small:
$$C(v,k,t) \drel{\geq} \frac{(t+1)(v-t)}{(k+1)(v-k)} {v \choose t}
    \bcdisplaydiv {k \choose t} \,.$$

The best general upper bound on $C(v,k,t)$ is due to R\"odl~\cite{rodl}:
Define the {\em density\/} of a covering to be the average number of blocks
containing a $t$-set.  The minimum density of a $(v,k,t)$ covering is
$C(v,k,t) {k \choose t} / {v \choose t}$ and is obviously at least~1. R{\"o}dl
shows that for $k$~and~$t$ fixed there exist coverings with density
approaching~1 as $v$~gets large. Erd{\H o}s and Spencer~\cite{erdos-spencer}
give the bound
\[ C(v,k,t) {k \choose t}  \bcdisplaydiv {v \choose t}
  \drel{\leq} 1 + \ln\!{k \choose t} \,,
\]
which is weaker but applies to all $v$,~$k$, and~$t$. Furthermore it can be
improved by at most a factor of $4 \ln 2 \approx 2.77$ asymptotically, because
a $(v, v\squish{-}1, \lfloor v/2 \rfloor)$ covering that achieves the
Sch{\"o}nheim lower bound has density asymptotic to~$v/4$, while the Erd{\H
o}s-Spencer upper bound in that case corresponds to a density asymptotic to $v
\ln 2$.

This paper presents new constructions for coverings. The greedy method of
Section~\ref{s:greedy} produces reasonably good coverings and it is
completely general---it applies to all possible values of $v$,~$k$, and~$t$,
and it doesn't rely on the existence of other good coverings. The finite
geometries of Section~\ref{s:geometries} produce very good (often optimal)
coverings, but they apply only to certain sets of $v$,~$k$, and~$t$ values.
The induced-covering method of Section~\ref{s:induced}, which constructs
coverings from larger ones, and the dynamic programming method of
Section~\ref{s:dynamic}, which constructs coverings from smaller ones, both
apply to all parameter values, but they rely on preexisting coverings. (We
show in a paper with Spencer~\cite{gkps} that the greedy construction, as
well as the induced-covering method applied to certain finite geometry
coverings, both produce coverings that match R{\"o}dl's bound.) Finally, the
previously known methods of Section~\ref{s:others}, when combined with the
methods of earlier sections, yield the tables of upper bounds in
Section~\ref{s:tables}.

\section{Greedy Coverings}
\label{s:greedy}

Our greedy algorithm for generating coverings is
analogous to the surprisingly good greedy algorithm of
Conway and Sloane~\cite{lex} for generating codes.
That algorithm may be stated very concisely: To construct a code
of length~$n$ and minimum distance~$d$, arrange the binary $n$-tuples
in lexicographic order, and repeatedly choose the first one in
the list that is distance $d$ or more from all $n$-tuples chosen earlier;
the $n$-tuples chosen are the codewords.  The resulting code is called
a {\em lexicographic code}, or {\em lexicode}.

This simple method has several nice features:
Lexicodes tend to be fairly good (at packing codewords into the space),
they are linear, and they include some well-known codes such
as Hamming codes and the binary Golay codes.
Brouwer, Shearer, Sloane, and Smith \cite[page~1349]{bsss} use the same method
to make constant weight codes, by choosing only $n$-tuples of a given weight.

The greedy algorithm does not require lexicographic order.
Brualdi and Pless~\cite{greedycodes} show that a large family of
orders lead to linear codes.
And sometimes Gray code orders, for example, lead to better codes.

Constructing good codes and good constant weight codes are packing problems.
But a similar method applies to covering problems.
A greedy $(v,k,t)$ covering is one generated by the following algorithm:
\begin{enumerate}
\item Arrange the $k$-subsets of a $v$-set in a list.

\item Choose from the list the $k$-subset that contains the maximum number of
$t$-sets that are still uncovered. In case of ties, choose the $k$-subset
occurring earliest in the list. \label{main-step}

\item Repeat Step~\ref{main-step} until all $t$-sets are covered.
\end{enumerate}

The list of $k$-sets can be in any order.
Some natural orders are lexicographic,
colex (which is similar to lexicographic
but the subsets are read from right to left rather than left to right),
and a generalized Gray code order
(where successive sets differ only by one deletion and one addition).
The resulting lists, when $k=3$ and $v=5$, are
\begin{small} \[ \begin{array}{l*{9}{@{\mspace{10}}l}l}
123 & 124 & 125 & 134 & 135 & 145 & 234 & 235 & 245 & 345 & \mbox{(lexicographic);} \\
123 & 124 & 134 & 234 & 125 & 135 & 235 & 145 & 245 & 345 & \mbox{(colex);} \\
123 & 134 & 234 & 124 & 145 & 245 & 345 & 135 & 235 & 125 & \mbox{(gray).}
\end{array} \] \end{small}
Nijenhuis and Wilf~\cite{nijenhuis-wilf} give
algorithms to generate lexicographic and Gray code orders.
Stanton and White~\cite{stanton-white} discuss colex algorithms.

It is natural to investigate the greedy algorithm with random order, too,
since we know~\cite{gkps} that random order does well asymptotically.
To keep with the constructive spirit of this paper,
we used an easily reproduced ``random'' permutation of the $k$-sets.
To generate the permutation,
start with the $k$-sets lexicographically ordered
in positions 1~through~${v \choose k}$,
then successively swap the $k$-sets in positions $i$~and $i+j$,
for $i=1$,~2, \dots,~${v \choose k}$,
where $j$~is $X_i \bmod \bigl( {v \choose k} - i + 1 \bigr)$
and where the sequence of pseudo-random~$X$'s
comes from the linear congruential generator
$X_{i+1} = (41 X_i + 7) \bmod 2^{30} \!.$
The seed~$X_0$ is~$1$,
and when there are multiple random-order runs
on the same set of $(v,k,t)$ parameters,
the subsequent seeds are 2,~3,~\dots.
Knuth~\cite{knuth-volume-2} discusses the linear congruential method.

Greedy coverings are not in general optimal, but as happens with codes
(Brouwer, Shearer, Sloane, and Smith~\cite{bsss}, Brualdi and
Pless~\cite{greedycodes}, Conway and Sloane~\cite{lex}) they are often quite
good---about 42\% of the table entries come from greedy coverings.
Interestingly, the Steiner system $S(24,8,5)$, which Conway and
Sloane~\cite[page~347]{lex} showed is a constant-weight lexicographic code,
also arises as a greedy covering.

The problem with greedy coverings is that they are expensive to compute.
Our implementation of the algorithm above uses two arrays:
one with $v \choose k$~locations corresponding to the $k$-subsets,
and one with $v \choose t$~locations corresponding to the $t$-subsets.
Each $k$-set array location contains the number of uncovered
$t$-sets contained in that $k$-set, and is initialized to~$k \choose t$.
Each $t$-set array location contains a $0$~or~$1$,
indicating whether that $t$-set has been covered.
Each time through Step~\ref{main-step},
each $t$-set contained in the selected $k$-set must be checked.
If the $t$-set is uncovered, it is marked as covered,
and each $k$-set containing it must have its array location decremented.
For fixed $k$~and~$t$, the algorithm
asymptotically takes time and space $O(v^k)$.

We ran a program to generate greedy coverings for
all entries in our tables, for all four orders described above.
For random order, we used $10^e$ runs,
where
$e = 3[v \relshrnk{\leq} 20] +
   [v \relshrnk{\leq} 15] + [v \relshrnk{\leq} 10] +
   [k \relshrnk{\leq} 10] + [k \relshrnk{\leq} 5] + 2[U]$
and where $U$~is the predicate `$t=2$ and $C(v,k,2)$ is unknown'
(the symbol $[P]$ is~1 if the predicate~$P$ is true, 0~otherwise).

For the range of parameters of our tables,
the four orders produced coverings of roughly the same size,
but lexicographic order performed slightly better on average than colex order,
which performed better than Gray code order,
which performed better than a single run of random order.

\section{Finite Geometry Coverings}
\label{s:geometries}

Finite geometries may be used to construct very good coverings for
certain sets of parameters.
Anderson~\cite{anderson} has a nice discussion of finite geometries.

Let $\PG(m,q)$ denote the projective geometry of dimension~$m$ over $\GF(q)$,
where $q$~is a prime power.
The points of $\PG(m,q)$ are the equivalence classes of
nonzero vectors $u=(u_0,u_1,\ldots,u_m)$, where two vectors $u$~and~$v$
are equivalent if $u = \lambda v$ for some nonzero $\lambda \in \GF(q)$.
There are $(q^{m+1}-1)/(q-1)$ such points.

A $k$-flat is a $k$-dimensional subspace of $\PG(m,q)$, for
$1\leq k \leq m$, determined by $m-k$ independent homogeneous linear
equations.
A $k$-flat has $(q^{k+1}-1)/(q-1)$ points,
and there are $\qbinomial{m+1}{k+1}$
different $k$-flats in $\PG(m,q)$,
where
$$
\qbinomial{n}{k} \deq
\frac{(q^n-1)(q^{n-1}-1) \ldots (q^{n-k+1}-1)}
{(q^k-1)(q^{k-1}-1) \ldots (q-1)}
$$
is the $q$-binomial coefficient.

By removing all points with $u_0 = 0$ we obtain the affine (or
Euclidean) geometry $\AG(m,q)$.
It has $q^m$~points and $q^{m-k} \qbinomial{m}{k}$
different $k$-flats, each of which contains $q^k$~points.

For either geometry,
any $k+1$ independent points determine a $k$-flat,
and $k+1$ dependent points are contained in multiple $k$-flats,
so the $k$-flats cover every set of $k+1$ points.
Thus,
taking the points of the geometry as the $v$-set of the covering,
and taking the points of a $k$-flat as a block of the covering,
we get the following two theorems.

\begin{theorem}\label{thm:proj}
$$
C\left(\frac{q^{m+1}-1}{q-1},\frac{q^{k+1}-1}{q-1},k\squish{+}1\right)
  \drel{\leq} \qbinomial{m+1}{k+1} \,.
$$
\end{theorem}

\begin{theorem}\label{thm:affine}
$$
C(q^m\!,q^k\!,k\squish{+}1)
  \drel{\leq} q^{m-k} \qbinomial{m}{k} \,.
$$
\end{theorem}

Equality holds for both theorems when $k = m-1$ or $k=1$.
Theorem~\ref{thm:proj} is due to Ray-Chaudhuri~\cite{ray},
and Theorem~\ref{thm:affine} follows easily from results of
Abraham, Ghosh, and Ray-Chaudhuri~\cite{abraham-et-al},
although the idea of using finite geometries to construct coverings
dates back at least to Veblen and Bussey~\cite{veblen-bussey} in 1906.

\section{Induced coverings}
\label{s:induced}

The main drawback of the finite geometry coverings is that
they exist only for certain families of parameters.
But they are such good coverings that they can be used to construct
pretty~good coverings for other parameters.

Suppose we have a good $(v,k,t)$ covering, say from a geometry, and
we want to construct a $(v',k',t)$ covering, where $v' < v$ and $k'<k$.
Consider the family of sets obtained from the ($k$-element) blocks
by randomly choosing $v'$~elements of the $v$-set,
deleting all other elements from the blocks,
and throwing out any blocks with fewer than $t$ elements
(since those blocks cover no $t$-sets).

The remaining blocks cover all $t$-subsets of the $v'$~elements, but
have different sizes.  Suppose some block has $\ell$~elements.
If $\ell=k'$ its size is correct as is,
and it becomes a block of our new covering.
If $\ell<k'$, add any $k'-\ell$ elements to the block.
And if $\ell>k'$, replace the block by an $(\ell,k',t)$ covering,
which covers all $t$-sets the original block covered.

The new blocks each have $k'$~elements, and together they cover all
$t$-sets, so the new family forms a $(v',k',t)$ {\em induced\/} covering.

In small cases, the method tends to do best when $k'/k$ is about $v'/v$. In
large cases, the method does well if for every $\ell$ near $v'k/v$, a good
$(\ell,k',t)$ covering is available. Also, it need not start with a finite
geometry covering---any $(v,k,t)$ covering will do. But generally the better
the covering it starts with, the better the result.

The induced coverings in our tables come either from using
the simple special cases of Section~\ref{s:trivial} or from finite geometries.
We constructed each finite geometry covering based on $PG(m,p)$ and $AG(m,p)$
with \mbox{$p \leq 11$} prime and with
at most $10^4$~points and $10^6$~flats.
For each such covering, and for each $v$~and~$k$ in the relevant table,
we used a random set of $v$~points
to construct an induced covering as described above,
trying 100 random sets in each case.

\section{Combining Smaller Coverings}
\label{s:dynamic}

Suppose we want to form a $(v_1\squish{+}v_2,k,t)$ covering. Let the
$(v_1{+}v_2)$-set be the disjoint union of a $v_1$-set and a $v_2$-set. Given
an~$s$ with $0 \leq s \leq t$, choose a $(v_1,\ell,s)$ covering and a
$(v_2,k\squish{-}\ell,t\squish{-}s)$ covering for some~$\ell$, which must be
in the range $s \leq \ell \leq k-t+s$. For each possible arrangement of
$t$~elements as an $s$-subset of the $v_1$-set and a $(t{-}s)$-subset of the
$v_2$-set, there is an $\ell$-set from the first covering and a
$(k{-}\ell)$-set from the second covering whose union is a $k$-set that
covers the $t$-set. Thus the number of blocks that cover all such $t$-sets is
at most the product of the sizes of the two coverings. Choosing an
optimal~$\ell$ for each~$s$ gives us our $(v_1\squish{+}v_2,k,t)$ covering
built up from smaller coverings. This construction gives the bound
\[ C(v_1\squish{+}v_2,k,t)
     \drel{\leq} \sum_{s=0}^t \min_\ell
                C(v_1,\ell,s) \cdot C(v_2,k\squish{-}\ell,t\squish{-}s) \,.
\]
Furthermore we can try all choices of $v_1$~and~$v_2$
summing to the $v$ of interest.

The coverings produced by this method tend to have some redundancy.
To remove redundancy when $v_1=2$, for example,
we can try combining a $(v,k,t)$ covering and
a $(2,0,0)$ covering (which has one block, the empty set),
along with a $(v,k\squish{-}2,t\squish{-}1)$ covering and a $(2,2,2)$ covering.
This forms a $(v\squish{+}2,k,t)$ covering,
and is sometimes an improvement over the basic construction above:
$$
C(v\squish{+}2,k,t)
     \drel{\leq} C(v,k,t) \,+\, C(v,k\squish{-}2,t\squish{-}1) \,.
$$

This example has replaced the $s$~and~$s+1$
terms of the basic construction's bound,
when $s=1$, with the single term
$$
\min_\ell C(v_1,\ell,s+1) \cdot C(v_2,k\squish{-}\ell,t\squish{-}s) \,.
$$
The new term corresponds to covering any $t$-subset having either
$s$~or~$s+1$ elements in the $v_1$-set, by using one product of coverings,
rather than two. If changing $C(v_1,\ell,s)$ to $C(v_1,\ell,s\squish{+}1)$
does not cost too much, the bound will improve.

To generalize this combining of terms, define $c_{i,j}$ for $0 \leq i \leq j
\leq t$ to be the number of blocks required to cover any $t$-subset that has
between $i$~and~$j$ elements in the $v_1$-set, and between $t-j$ and $t-i$
elements in the $v_2$-set.  Since $c_{i,j} \leq c_{i,r} + c_{r+1,j}$ for any
$i \leq r < j$, we have
\begin{multline*}
c_{i,j} \drel{\leq} \min \bigl(\min_\ell C(v_1,\ell,j) \cdot C(v_2,k\squish{-}\ell,t\squish{-}i), \\
    \min_{i \leq r < j} (c_{i,r} \squish{+} c_{r+1,j}) \bigr) \,.
\end{multline*}
Using dynamic programming, we may efficiently compute a bound for~$c_{0,t}$,
which is an upper bound for $C(v_1\squish{+}v_2,k,t)$.

This general construction produces about 30\% of the entries in our tables.
It includes as special cases several of the simple constructions of
Section~\ref{s:trivial}, as well as the direct-product construction of Morley
and van Rees~\cite{mr}, which yields the bound
\[ C(2v\squish{+}y,v\squish{+}k\squish{+}y,t\squish{+}s\squish{+}1)
   \drel{\leq} C(v,k,t) \,+\, C(v\squish{+}y,k\squish{+}y,s) \,.
\]

\section{Other Constructions}
\label{s:others}

\subsection{Simple Constructions}
\label{s:trivial}

There are several simple and well-known methods
for building coverings from other coverings.
All but the last of these methods are special cases of
the methods in the previous two sections.

Adding a random element to each block of a $(v,k,t)$ covering
gives a $(v,k\squish{+}1,t)$ covering of the same size.  Thus
\[ C(v,k\squish{+}1,t) \drel{\leq} C(v,k,t) \,.
\]
Adding a new element to a $v$-set,
and including it in every block in a $(v,k,t)$ covering,
forms a $(v\squish{+}1,k\squish{+}1,t)$ covering of the same size, hence
\[ C(v\squish{+}1,k\squish{+}1,t) \drel{\leq} C(v,k,t) \,.
\]
Combining a $(v,k,t)$ covering and a $(v,k\squish{-}1,t\squish{-}1)$ covering
over the same $v$-set, by adding a new $v{+}1$st element to
all of the blocks of the $(v,k\squish{-}1,t\squish{-}1)$ covering
but to none of the blocks of the $(v,k,t)$ covering,
forms a $(v\squish{+}1,k,t)$ covering,
of size the sum of the other two sizes, thus
\[ C(v\squish{+}1,k,t)
    \drel{\leq} C(v,k,t) \,+\, C(v,k\squish{-}1,t\squish{-}1) \,.
\]
Those constructions are special cases of
the method of Section~\ref{s:dynamic}.

Deleting one element from a $v$-set, and adding a random element to any
block of a $(v,k,t)$ covering that contains the deleted element,
creates a $(v\squish{-}1,k,t)$ covering of the same size.  Thus
\[ C(v\squish{-}1,k,t) \drel{\leq} C(v,k,t) \,.
\]
Choosing the element of a covering that occurs in the fewest blocks,
throwing away all other blocks, and then throwing away the chosen element,
results in a $(v\squish{-}1,k\squish{-}1,t\squish{-}1)$ covering.
This method, due to Sch{\"o}nheim,
is a reformulation of Theorem~\ref{schonheim-one-level};
the corresponding upper bound is
\[ C(v\squish{-}1,k\squish{-}1,t\squish{-}1)
        \drel{\leq} \Bigl\lfloor \frac{k}{v} \, C(v,k,t) \Bigr\rfloor \,.
\]
Those two constructions are special cases of the induced-covering method
of Section~\ref{s:induced}.

Replacing each element of the $v$-set in a $(v,k,t)$ covering
by $m$~different elements
gives an $(mv,mk,t)$ covering of the same size, thus
\[ C(mv,mk,t) \drel{\leq} C(v,k,t) \,.
\]

\subsection{Steiner Systems}
\label{s:steiner}

A {\em Steiner system\/} is a covering in which the covering density
is~1---every $t$-set is covered exactly once.  Clearly a Steiner system is an
optimal covering, as well as an optimal packing, and $C(v,k,t)=L(v,k,t)$. The
projective and affine coverings by lines (1-flats), for example, are Steiner
systems. Brouwer, Shearer, Sloane, and Smith~\cite[page~1342]{bsss} and Chee,
Colbourn, and Kreher~\cite{chee-colbourn-kreher} give tables of small Steiner
systems.

If a $(v,k,t)$ Steiner system exists then
$C(v\squish{+}1,k,t)=L(v\squish{+}1,k,t)$. This result is due to
Sch\"onheim~\cite[Theorem~II]{schonheim}; the proof also appears in Mills and
Mullin~\cite[Theorem~1.3]{mills-mullin}.

\subsection{Tur\'an Theory}
\label{s:turan}

The {\em Tur\'an number\/} $T(n,\ell,r)$ is the minimum number of $r$-subsets
of an $n$-set such that every $\ell$-subset
contains at least one of the $r$-subsets.
It is easy to see that
\[ C(v,k,t) \deq T(v,v\squish{-}t,v\squish{-}k) \,,
\]
so covering numbers are just Tur\'an numbers reordered. The two sets of
numbers, however, have been studied for different parameter ranges (de~Caen's
lower bound in the introduction, for instance, is useful primarily for
Tur\'an theory ranges). Most papers on coverings have $v$~large compared with
$k$~and~$t$, while most papers on Tur\'an numbers have $n$~large compared
with $\ell$~and~$r$, often focusing on the quantity $\lim_{n \rightarrow
\infty} T(n,\ell,r) / {n \choose r}$ for fixed $\ell$~and~$r$. Thus Tur\'an
theory usually studies $C(v,k,t)$ for $k$~and~$t$ not too far from~$v$.

Fifty years ago Tur\'an~\cite{turan} determined $T(n,\ell,2)$ exactly,
showing that $C(v,v\squish{-}2,t) = L(v,v\squish{-}2,t)$, the Sch\"onheim
lower bound. He also gave upper bounds and conjectures for $T(n,4,3)$ and
$T(n,5,3)$, which stimulated much of the research. The results labeled
`Tur\'an theory' in our tables either are described in recent survey papers
by de~Caen~\cite{decaen} and Sidorenko~\cite{sidorenko}, or follow from
constructions due to de~Caen, Kreher, and
Wiseman~\cite{decaen-kreher-wiseman} \nocite{decaen-krm} or to
Sidorenko~\cite{sidorenko-email}.

Sidorenko~\cite{sidorenko-email} also recently told us of a Tur\'an theory
construction, similar in spirit to the combining constructions of
Section~\ref{s:dynamic}, that improves many bounds in the table. In terms of
covering theory, let $x$ be an element occurring in the most blocks of a
$(v,k,t)$ covering, and replace~$x$ by $x'$~and~$x''$: If a block~$b$ did not
contain~$x$, replace it by two blocks, $b \cup \{x'\}$ and $b \cup \{x''\}$;
if $b$ did contain~$x$, replace it by the single block $b - \{x\} \cup
\{x',x''\}$. Finally, add a $(v\squish{-}1,k\squish{+}1,t\squish{+}1)$
covering on the same elements minus $x'$~and~$x''\!$. It is not hard to see
that this is a $(v\squish{+}1,k\squish{+}1,t\squish{+}1)$ covering, and that
it gives the bound
\begin{multline*}
C(v\squish{+}1,k\squish{+}1,t\squish{+}1) \\
     \drel{\leq} \lfloor (2v-k) \, C(v,k,t) / v \rfloor
     \,+\, C(v\squish{-}1,k\squish{+}1,t\squish{+}1) \,.
\end{multline*}

\subsection{Cyclic Coverings}

Another well-known method that is often successful when applicable---when the
size of a prospective covering is~$v$---is to construct a cyclic covering:
Choose some $k$-subset as the first block, and choose the $v-1$ cyclic shifts
of that block as the remaining blocks. Trying this for all possible $k$-sets
is fairly cheap, and frequently it produces a covering. The entries
$C(19,9,3) \leq 19$ and $C(24,10,3) = 24$ in our tables, for example, are
generated by the $k$-sets 1~2~3~4~6~8~13~14~17 and 1~2~3~5~6~8~12~13~15~21,
and are unmatched by any other method.

Incidentally, if the size of a prospective covering is a multiple of~$v$,
say~$2v$, the same method applies by taking
the cyclic shifts of two starting blocks;
the few cases we tried for this variation
produced no improvements in the tables.

\subsection{Hill-Climbing}
\label{s:hill-climbing}

For cases of interest---with $v$ not too large---random coverings are not
very good, but hill-climbing sometimes finds good coverings: Start with a
fixed number of random $k$-sets, say $L(v,k,t)+\epsilon$ for some small
integer~$\epsilon$. Rank the $k$-sets by the number of $t$-sets they cover
that no other $k$-set covers, and replace one with lowest rank by another
random $k$-set. Repeat until all $t$-sets are covered or until time runs out.

We found a few good coverings with this method, but Nurmela and
\"{O}sterg{\aa}rd~\cite{no} went much further, using simulated annealing---a
more sophisticated hill-climbing---to find many good coverings. In fact many
of the bounds in the tables could be improved, by starting with a covering
produced by one of the other methods and then hill-climbing; but generally
the improvements would be small.

\section{Tables of Upper Bounds on $C(v,k,t)$}
\label{s:tables}

We constructed Tables \ref{t=2} through~\ref{t=8}
using the methods described above,
together with results from the literature.
Each table entry indicates the upper bound, the method of construction, and
whether the covering is known to be optimal.
We have tried to provide constructions for as many sets of parameters
as possible, so we list a method of construction from this paper
even when a result in the literature achieves the same bound.
When two different methods produce the same size covering,
we've given precedence to the method listed earlier in the Key to the tables.

About 93\% of the 1631 nontrivial ($v \relshrnk{>} k \relshrnk{>} t$)
upper bounds in the tables
come from one of the constructions described in this paper.
For each of the remaining upper bounds,
there is a source in our reference list that describes the result,
although to keep our reference list reasonably short
we have often given a secondary source rather than the original.
(Mills and Mullin~\cite{mills-mullin} give an extensive list
of previous results and references.)
Sources for Steiner systems, Tur\'an number bounds,
and simulated annealing coverings
appear in Sections \ref{s:steiner},~\ref{s:turan}, and~\ref{s:hill-climbing};
the Todorov constructions come from papers by
Todorov~\cite{todorov-thesis,todorov-13-blocks,todorov-pairs} and
Todorov and Tonchev~\cite{todorov-tonchev};
and the remaining upper bounds appear in Table~\ref{tab:misc}.
The covering number $C(24,18,17)$ is listed in Table~\ref{tab:misc},
even though it doesn't occur in the other tables,
because it yields a $(15,9,8)$ simple induced covering
(of Section~\ref{s:trivial}).

Gordon et~al.~\cite{gppt} construct an optimal $(12,6,3)$ covering,
using a block-array construction.
That method directly extends to the $(18,9,4)$ covering given in
Table~\ref{tab:misc}, and a similar construction
gives four other coverings listed in the table.

\begin{table}
\begin{center}
\begin{tabular}{|c|l|}
\hline
bound & reference \\ \hline
$C(29,5,2) \leq 44 $ & Lamken \cite{lmmv}\\
$C(31,7,2) = 26 $ & Todorov~\cite{todorov-pairs} techniques (lower bound)%
\nocite{todorov-t=2-lower-bounds}\\
$C(12,6,3) = 15 $ & Gordon et al.\ \cite{gppt}\\
$C(14,6,3) \leq 25$ & Lotto covering \cite{lotto}\\
$C(15,6,3) \leq 31$ & Lotto covering \cite{lotto}\\
$C(16,6,3) \leq 38 $ & Hoehn \cite{hoehn}\\
$C(18,6,3) = 48$ & Lotto covering \cite{lotto}\\
$C(30,6,3) \leq 237$ & Lotto covering \cite{lotto}\\
$C(11,7,4) = 17$ & Sidorenko \cite{sidorenko-email}\\
$C(14,6,4) \leq 87$ & Hoehn \cite{hoehn}\\
$C(18,6,4) \leq 258$ & Lotto covering \cite{lotto}\\
$C(18,9,4) \leq 43 $ & Gordon et al.\ \cite{gppt}\\
$C(20,10,4) \leq 43$ & block-array construction\\
$C(24,12,5) \leq 86$ & block-array construction\\
$C(30,15,5) \leq 120$ & block-array construction\\
$C(12,8,6) \leq 51 $ & Morley \cite{morley}\\
$C(32,16,6) \leq 286$ & block-array construction\\
$C(15,12,8) = 30 $ & Radziszowski and Sidorenko \cite{rad-sid} \\
$C(24,18,17) = 21252 $ & de Caen \cite{decaen} \\
\hline
\end{tabular}
\end{center}
\caption{Miscellaneous results}\label{tab:misc}
\end{table}

Most of the lower bounds used to establish optimality
follow from the Sch{\"o}nheim inequality (Theorem~\ref{schonheim-one-level});
and a few others are listed as equalities in Table~\ref{tab:misc}.
For the rest:
If $t=2$, the lower bound is explained by Mills and Mullin~\cite{mills-mullin}
when it is less than~14 or has $v \leq 5$,
or explained by Todorov~\cite{todorov-pairs} otherwise;
if $t=3$, it's either Mills and Mullin or
Todorov and Tonchev~\cite{todorov-tonchev};
and if $4 \leq t \leq 8$,
it's either Mills~\cite[Theorem~2.3]{mills-lower-bounds},
Todorov~\cite[Theorem~4]{todorov-lower-bounds}, or
Sidorenko's Tur\'an theory survey~\cite{sidorenko}.

How good are our bounds? For $t=2$, very good---most of the entries are known
to be optimal, and the largest gap between an entry's lower and upper bound
is currently only a factor of 1.12. That largest gap rises with~$t$, though,
to 1.89 for $t=4$, to 2.98 for $t=6$, and to $3.72$ for $t=8$. We believe
that our lower bounds tend to be closer to the truth than our upper bounds;
it's quite possible that all the upper bounds are within a factor of~3, but
probably not a factor of~2, of optimal.

Most of the entries in the tables for $t>2$ are not optimal,
and we would appreciate knowing of any better coverings.
Please send communications to the first author, at
{\tt gordon@ccrwest.org}.

\bigskip
\medskip

\noindent
Key to Tables \ref{t=2} through \ref{t=8}
\nopagebreak

\vspace{1ex}
\noindent
\begin{tabular}{c@{$\mspace{4}$---$\mspace{6}$}l@{}}
$l$ & greedy covering, lexicographic order\\ 
$c$ & greedy covering, colex order\\ 
$g$ & greedy covering, Gray code order\\ 
$r$ & greedy covering, random order\\ 
$p$ & projective geometry covering\\ 
$a$ & affine geometry covering\\ 
$o$ & cyclic covering\\ 
$m$ & multiple of smaller covering\\ 
$e$ & simple dynamic programming  (Section~\ref{s:trivial})\\ 
$j$ & simple induced covering (Section~\ref{s:trivial})\\ 
$d$ & dynamic programming method (Section~\ref{s:dynamic})\\ 
$i$ & induced covering\\ 
$u$ & Sidorenko Tur\'an construction (Section~\ref{s:turan})\\ 
$s$ & Steiner system\\ 
$t$ & Tur\'an theory\\ 
$x$ & covering with small $k$ and $t$; see Mills and Mullin \cite[\S 3]{mills-mullin}\\ 
$y$ & covering with fixed size; see Mills and Mullin \cite[\S 4]{mills-mullin}\\ 
$v$ & Todorov construction\\ 
$w$ & was known previously; see Table \ref{tab:misc}\\ 
$n$ & Nurmela-\"{O}sterg{\aa}rd simulated annealing covering\\ 
$h$ & hill-climbing\\ 
$*$ & optimal covering\\ 
\end{tabular}
\begin{table}
\begin{tabular}{r|rrrrrrrrrrrrrr}
$v \backslash k$ & 3 & 4 & 5 & 6 & 7 & 8 & 9 & 10 & 11 & 12 & 13 & 14 & 15 & 16 \\ \hline 
3 & ${1}$\rlap{$\null^{ *}$} &  & & & & & & & & & & & &\\ 
4 & ${3}$\rlap{$\null^{l*}$} & ${1}$\rlap{$\null^{ *}$} &  & & & & & & & & & & &\\ 
5 & ${4}$\rlap{$\null^{l*}$} & ${3}$\rlap{$\null^{l*}$} & ${1}$\rlap{$\null^{ *}$} &  & & & & & & & & & &\\ 
6 & ${6}$\rlap{$\null^{o*}$} & ${3}$\rlap{$\null^{l*}$} & ${3}$\rlap{$\null^{l*}$} & ${1}$\rlap{$\null^{ *}$} &  & & & & & & & & &\\ 
7 & ${7}$\rlap{$\null^{l*}$} & ${5}$\rlap{$\null^{l*}$} & ${3}$\rlap{$\null^{l*}$} & ${3}$\rlap{$\null^{l*}$} & ${1}$\rlap{$\null^{ *}$} &  & & & & & & & &\\ 
8 & ${11}$\rlap{$\null^{l*}$} & ${6}$\rlap{$\null^{l*}$} & ${4}$\rlap{$\null^{l*}$} & ${3}$\rlap{$\null^{l*}$} & ${3}$\rlap{$\null^{l*}$} & ${1}$\rlap{$\null^{ *}$} &  & & & & & & &\\ 
9 & ${12}$\rlap{$\null^{r*}$} & ${8}$\rlap{$\null^{l*}$} & ${5}$\rlap{$\null^{l*}$} & ${3}$\rlap{$\null^{l*}$} & ${3}$\rlap{$\null^{l*}$} & ${3}$\rlap{$\null^{l*}$} & ${1}$\rlap{$\null^{ *}$} &  & & & & & &\\
\end{tabular}
\caption{$t=2$}

\label{t=2}

\end{table}
\begin{table}
\begin{tabular}{r|rrrrrrrrrrrrr}
$v \backslash k$ & 4 & 5 & 6 & 7 & 8 & 9 & 10 & 11 & 12 & 13 & 14 & 15 & 16 \\ \hline 
4 & ${1}$\rlap{$\null^{ *}$} &  & & & & & & & & & & &\\ 
5 & ${4}$\rlap{$\null^{l*}$} & ${1}$\rlap{$\null^{ *}$} &  & & & & & & & & & &\\ 
6 & ${6}$\rlap{$\null^{o*}$} & ${4}$\rlap{$\null^{l*}$} & ${1}$\rlap{$\null^{ *}$} &  & & & & & & & & &\\ 
7 & ${12}$\rlap{$\null^{r*}$} & ${5}$\rlap{$\null^{l*}$} & ${4}$\rlap{$\null^{l*}$} & ${1}$\rlap{$\null^{ *}$} &  & & & & & & & &\\ 
8 & ${14}$\rlap{$\null^{l*}$} & ${8}$\rlap{$\null^{o*}$} & ${4}$\rlap{$\null^{l*}$} & ${4}$\rlap{$\null^{l*}$} & ${1}$\rlap{$\null^{ *}$} &  & & & & & & &\\ 
9 & ${25}$\rlap{$\null^{l*}$} & ${12}$\rlap{$\null^{l*}$} & ${7}$\rlap{$\null^{j*}$} & ${4}$\rlap{$\null^{l*}$} & ${4}$\rlap{$\null^{l*}$} & ${1}$\rlap{$\null^{ *}$} &  & & & & & &\\ 
10 & ${30}$\rlap{$\null^{r*}$} & ${17}$\rlap{$\null^{r*}$} & ${10}$\rlap{$\null^{l*}$} & ${6}$\rlap{$\null^{l*}$} & ${4}$\rlap{$\null^{l*}$} & ${4}$\rlap{$\null^{l*}$} & ${1}$\rlap{$\null^{ *}$} &  & & & & &\\
\end{tabular}
\caption{$t=3$}

\end{table}
\begin{table}
\begin{tabular}{r|rrrrrrrrrrrr}
$v \backslash k$ & 5 & 6 & 7 & 8 & 9 & 10 & 11 & 12 & 13 & 14 & 15 & 16 \\ \hline 
5 & ${1}$\rlap{$\null^{ *}$} &  & & & & & & & & & &\\ 
6 & ${5}$\rlap{$\null^{l*}$} & ${1}$\rlap{$\null^{ *}$} &  & & & & & & & & &\\ 
7 & ${9}$\rlap{$\null^{l*}$} & ${5}$\rlap{$\null^{l*}$} & ${1}$\rlap{$\null^{ *}$} &  & & & & & & & &\\ 
8 & ${20}$\rlap{$\null^{r*}$} & ${7}$\rlap{$\null^{j*}$} & ${5}$\rlap{$\null^{l*}$} & ${1}$\rlap{$\null^{ *}$} &  & & & & & & &\\ 
9 & ${30}$\rlap{$\null^{r*}$} & ${12}$\rlap{$\null^{l*}$} & ${6}$\rlap{$\null^{l*}$} & ${5}$\rlap{$\null^{l*}$} & ${1}$\rlap{$\null^{ *}$} &  & & & & & &\\ 
10 & ${51}$\rlap{$\null^{r}$} & ${20}$\rlap{$\null^{r*}$} & ${10}$\rlap{$\null^{o*}$} & ${5}$\rlap{$\null^{l*}$} & ${5}$\rlap{$\null^{l*}$} & ${1}$\rlap{$\null^{ *}$} &  & & & & &\\ 
11 & ${66}$\rlap{$\null^{j*}$} & ${32}$\rlap{$\null^{n*}$} & ${17}$\rlap{$\null^{j*}$} & ${9}$\rlap{$\null^{j*}$} & ${5}$\rlap{$\null^{l*}$} & ${5}$\rlap{$\null^{l*}$} & ${1}$\rlap{$\null^{ *}$} &  & & & &\\
\end{tabular}
\caption{$t=4$}

\end{table}
\begin{table}
\begin{tabular}{r|rrrrrrrrrrr}
$v \backslash k$ & 6 & 7 & 8 & 9 & 10 & 11 & 12 & 13 & 14 & 15 & 16 \\ \hline 
6 & ${1}$\rlap{$\null^{ *}$} &  & & & & & & & & &\\ 
7 & ${6}$\rlap{$\null^{l*}$} & ${1}$\rlap{$\null^{ *}$} &  & & & & & & & &\\ 
8 & ${12}$\rlap{$\null^{l*}$} & ${6}$\rlap{$\null^{l*}$} & ${1}$\rlap{$\null^{ *}$} &  & & & & & & &\\ 
9 & ${30}$\rlap{$\null^{r*}$} & ${9}$\rlap{$\null^{o*}$} & ${6}$\rlap{$\null^{l*}$} & ${1}$\rlap{$\null^{ *}$} &  & & & & & &\\ 
10 & ${50}$\rlap{$\null^{r*}$} & ${20}$\rlap{$\null^{d*}$} & ${8}$\rlap{$\null^{j*}$} & ${6}$\rlap{$\null^{l*}$} & ${1}$\rlap{$\null^{ *}$} &  & & & & &\\ 
11 & ${100}$\rlap{$\null^{n}$} & ${34}$\rlap{$\null^{j}$} & ${16}$\rlap{$\null^{j*}$} & ${7}$\rlap{$\null^{l*}$} & ${6}$\rlap{$\null^{l*}$} & ${1}$\rlap{$\null^{ *}$} &  & & & &\\ 
12 & ${132}$\rlap{$\null^{s*}$} & ${59}$\rlap{$\null^{n}$} & ${26}$\rlap{$\null^{t*}$} & ${12}$\rlap{$\null^{o*}$} & ${6}$\rlap{$\null^{l*}$} & ${6}$\rlap{$\null^{l*}$} & ${1}$\rlap{$\null^{ *}$} &  & & &\\
\end{tabular}
\caption{$t=5$}

\end{table}
\begin{table}
\begin{tabular}{r|rrrrrrrrrr}
$v \backslash k$ & 7 & 8 & 9 & 10 & 11 & 12 & 13 & 14 & 15 & 16 \\ \hline 
7 & ${1}$\rlap{$\null^{ *}$} &  & & & & & & & &\\ 
8 & ${7}$\rlap{$\null^{l*}$} & ${1}$\rlap{$\null^{ *}$} &  & & & & & & &\\ 
9 & ${16}$\rlap{$\null^{l*}$} & ${7}$\rlap{$\null^{l*}$} & ${1}$\rlap{$\null^{ *}$} &  & & & & & &\\ 
10 & ${45}$\rlap{$\null^{r*}$} & ${12}$\rlap{$\null^{j*}$} & ${7}$\rlap{$\null^{l*}$} & ${1}$\rlap{$\null^{ *}$} &  & & & & &\\ 
11 & ${84}$\rlap{$\null^{j}$} & ${29}$\rlap{$\null^{t*}$} & ${10}$\rlap{$\null^{j*}$} & ${7}$\rlap{$\null^{l*}$} & ${1}$\rlap{$\null^{ *}$} &  & & & &\\ 
12 & ${177}$\rlap{$\null^{n}$} & ${51}$\rlap{$\null^{w}$} & ${22}$\rlap{$\null^{d*}$} & ${9}$\rlap{$\null^{j*}$} & ${7}$\rlap{$\null^{l*}$} & ${1}$\rlap{$\null^{ *}$} &  & & &\\ 
13 & ${264}$\rlap{$\null^{n}$} & ${104}$\rlap{$\null^{n}$} & ${40}$\rlap{$\null^{t}$} & ${16}$\rlap{$\null^{d*}$} & ${8}$\rlap{$\null^{l*}$} & ${7}$\rlap{$\null^{l*}$} & ${1}$\rlap{$\null^{ *}$} &  & &\\
\end{tabular}
\caption{$t=6$}

\end{table}
\begin{table}
\begin{tabular}{r|rrrrrrrrr}
$v \backslash k$ & 8 & 9 & 10 & 11 & 12 & 13 & 14 & 15 & 16 \\ \hline 
8 & ${1}$\rlap{$\null^{ *}$} &  & & & & & & &\\ 
9 & ${8}$\rlap{$\null^{l*}$} & ${1}$\rlap{$\null^{ *}$} &  & & & & & &\\ 
10 & ${20}$\rlap{$\null^{j*}$} & ${8}$\rlap{$\null^{l*}$} & ${1}$\rlap{$\null^{ *}$} &  & & & & &\\ 
11 & ${63}$\rlap{$\null^{j*}$} & ${15}$\rlap{$\null^{l*}$} & ${8}$\rlap{$\null^{l*}$} & ${1}$\rlap{$\null^{ *}$} &  & & & &\\ 
12 & ${126}$\rlap{$\null^{t}$} & ${40}$\rlap{$\null^{r*}$} & ${12}$\rlap{$\null^{o*}$} & ${8}$\rlap{$\null^{l*}$} & ${1}$\rlap{$\null^{ *}$} &  & & &\\ 
13 & ${297}$\rlap{$\null^{n}$} & ${79}$\rlap{$\null^{n}$} & ${30}$\rlap{$\null^{d}$} & ${11}$\rlap{$\null^{j*}$} & ${8}$\rlap{$\null^{l*}$} & ${1}$\rlap{$\null^{ *}$} &  & &\\ 
14 & ${474}$\rlap{$\null^{j}$} & ${183}$\rlap{$\null^{e}$} & ${58}$\rlap{$\null^{t}$} & ${22}$\rlap{$\null^{d*}$} & ${10}$\rlap{$\null^{j*}$} & ${8}$\rlap{$\null^{l*}$} & ${1}$\rlap{$\null^{ *}$} &  &\\
\end{tabular}
\caption{$t=7$}

\end{table}
\begin{table}
\begin{tabular}{r|rrrrrrrr}
$v \backslash k$ & 9 & 10 & 11 & 12 & 13 & 14 & 15 & 16 \\ \hline 
9 & ${1}$\rlap{$\null^{ *}$} &  & & & & & &\\ 
10 & ${9}$\rlap{$\null^{l*}$} & ${1}$\rlap{$\null^{ *}$} &  & & & & &\\ 
11 & ${25}$\rlap{$\null^{r*}$} & ${9}$\rlap{$\null^{l*}$} & ${1}$\rlap{$\null^{ *}$} &  & & & &\\ 
12 & ${84}$\rlap{$\null^{t*}$} & ${18}$\rlap{$\null^{l*}$} & ${9}$\rlap{$\null^{l*}$} & ${1}$\rlap{$\null^{ *}$} &  & & &\\ 
13 & ${185}$\rlap{$\null^{t}$} & ${52}$\rlap{$\null^{t*}$} & ${15}$\rlap{$\null^{d*}$} & ${9}$\rlap{$\null^{l*}$} & ${1}$\rlap{$\null^{ *}$} &  & &\\ 
14 & ${482}$\rlap{$\null^{e}$} & ${121}$\rlap{$\null^{u}$} & ${40}$\rlap{$\null^{d}$} & ${13}$\rlap{$\null^{d*}$} & ${9}$\rlap{$\null^{l*}$} & ${1}$\rlap{$\null^{ *}$} &  &\\ 
15 & ${790}$\rlap{$\null^{j}$} & ${300}$\rlap{$\null^{u}$} & ${81}$\rlap{$\null^{t}$} & ${30}$\rlap{$\null^{d*}$} & ${12}$\rlap{$\null^{d*}$} & ${9}$\rlap{$\null^{l*}$} & ${1}$\rlap{$\null^{ *}$} & \\
\end{tabular}
\caption{$t=8$}

\label{t=8}

\end{table}

\clearpage

\acknowledgments

We thank W.~H. Mills, Nick Patterson, Alexander Sidorenko, D.~T. Todorov,
and an anonymous referee for some constructive suggestions
and for pointing out literature results we were unaware of.
We are particularly grateful to Alexander Sidorenko for
allowing us to publish his Tur\'an construction in Section~\ref{s:turan}.


\newcommand{\noopsort}[1]{} \newcommand{\eatperiod}[1]{}

\end{document}